\documentclass[openany, amssymb, psamsfonts]{amsart}
\usepackage{mathrsfs,comment}
\usepackage{xcolor}
\usepackage[normalem]{ulem}
\usepackage{url}
\usepackage[all,arc,2cell]{xy}
\UseAllTwocells
\usepackage{enumerate}
\usepackage{esint}  
\usepackage{hyperref}  
\usepackage{tikz}
\hypersetup{%
  bookmarksnumbered=true,%
  bookmarks=true,%
  colorlinks=true,%
  linkcolor=blue,%
  citecolor=blue,%
  filecolor=blue,%
  menucolor=blue,%
  pagecolor=blue,%
  urlcolor=blue,%
  pdfnewwindow=true,%
  pdfstartview=FitBH}   
\theoremstyle{definition}

\newtheorem{remark}{Remark}

\makeatletter
\let\c@obs=\c@thm
\let\c@cor=\c@thm
\let\c@prop=\c@thm
\let\c@lem=\c@thm
\let\c@prob=\c@thm
\let\c@con=\c@thm
\let\c@conj=\c@thm
\let\c@defn=\c@thm
\let\c@notn=\c@thm
\let\c@notns=\c@thm
\let\c@exmp=\c@thm
\let\c@ax=\c@thm
\let\c@pro=\c@thm
\let\c@ass=\c@thm
\let\c@warn=\c@thm
\let\c@rem=\c@thm
\let\c@sch=\c@thm
\let\c@equation\c@thm
\numberwithin{equation}{section}
\makeatother

\title{Solving the Laplace Equation and Applications in Imaging}
\author{Arina Oberoi}
\date{December 2024}

\begin{document}

\maketitle

\begin{abstract}
This paper examines solutions to the Laplace equation using analytical techniques, including separation of variables and the Poisson integral formula, and probabilistic methods, such as Brownian motion. We address applications to imaging, including noise reduction, data interpolation, and resolution enhancement, as well as discuss theoretical connections between the Laplace equation, specifically concepts such as harmonic functions, the mean value property, and the maximum principle, to provide a mathematical framework for practical implementation in imaging workflows.
\end{abstract}

\tableofcontents

\section{Solving the Laplace Equation in a Disc} 

The Laplace equation in Cartesian coordinates is 
\[
\nabla^2 u = \frac{\partial^2 u}{\partial x^2} + \frac{\partial^2 u}{\partial y^2} = 0
\]
where \(x = r\cos(\theta)\) and \(y = r\sin(\theta)\). To solve the Laplace equation \(\nabla^2 u = 0\) in a disc, we convert the equation to polar coordinates \((r,\theta)\)
\[
\frac{1}{r} \frac{\partial}{\partial r} \left( r \frac{\partial u}{\partial r} \right) + \frac{1}{r^2} \frac{\partial^2 u}{\partial \theta^2} = 0
\]
Assume the boundary condition \(u(R,\theta) = f(\theta)\) on the disc of radius \(R\) and apply separation of variables by writing \(u(r,\theta) = R(r)\Theta(\theta)\), which splits the equation into radial and angular parts. Substituting \(u(r,\theta) = R(r)\Theta(\theta)\) into the Laplace equation, we have
\[
\Theta(\theta) \left[ \frac{1}{r} \frac{d}{dr} \left( r \frac{dR}{dr} \right) \right] + R(r) \frac{1}{r^2} \left[ \frac{d^2 \Theta}{d\theta^2} \right] = 0
\]
\[
\frac{r^2}{R(r)} \left[ \frac{1}{r} \frac{d}{dr} \left( r \frac{dR}{dr} \right) \right] + \frac{1}{\Theta(\theta)} \left[ \frac{d^2 \Theta}{d\theta^2} \right] = 0
\]
\[
\frac{r}{R(r)} \frac{d}{dr} \left( r \frac{dR}{dr} \right) + \frac{1}{\Theta(\theta)} \frac{d^2 \Theta}{d\theta^2} = 0
\]
Since the first term depends only on \(r\) and the second term depends only on \(\theta\), and \(r\) and \(\theta\) are independent variables, each term must be equal to a constant that sums to zero. Therefore, we set
\[
\frac{r}{R(r)} \frac{d}{dr} \left( r \frac{dR}{dr} \right) = -\lambda, \quad \frac{1}{\Theta(\theta)} \frac{d^2 \Theta}{d\theta^2} = \lambda\]
where \(\lambda\) is the separation constant. To determine the sign and possible values of \(\lambda\), we consider the angular equation
\[
\frac{d^2 \Theta}{d\theta^2} + \lambda \Theta = 0
\]
\[
\int_0^{2\pi} \Theta \frac{d^2 \Theta}{d\theta^2} \, d\theta + \lambda \int_0^{2\pi} \Theta^2 \, d\theta = 0
\]
\[
\int_0^{2\pi} \Theta \frac{d^2 \Theta}{d\theta^2} \, d\theta = \left[ \Theta \frac{d\Theta}{d\theta} \right]_0^{2\pi} - \int_0^{2\pi} \left( \frac{d\Theta}{d\theta} \right)^2 \, d\theta
\]
For periodic boundary conditions, \(\Theta(0) = \Theta(2\pi)\) and \(\frac{d\Theta}{d\theta}(0) = \frac{d\Theta}{d\theta}(2\pi)\), so the boundary terms disappear
\[
\left[ \Theta \frac{d\Theta}{d\theta} \right]_0^{2\pi} = 0
\]
\[
-\int_0^{2\pi} \left( \frac{d\Theta}{d\theta} \right)^2 \, d\theta + \lambda \int_0^{2\pi} \Theta^2 \, d\theta = 0
\]
\[
\lambda \int_0^{2\pi} \Theta^2 \, d\theta = \int_0^{2\pi} \left( \frac{d\Theta}{d\theta} \right)^2 \, d\theta
\]
Since both \(\Theta^2\) and \(\left( \frac{d\Theta}{d\theta} \right)^2\) are non-negative, it follows that \(\lambda \geq 0\) Because the periodic boundary conditions require that \(\Theta(\theta)\) be sinusoidal functions (e.g., \(\cos(n\theta)\) or \(\sin(n\theta)\)) where \(n\) is an integer, we have that \(\lambda\) is be a non-negative integer.

\subsubsection*{Does a classical \(C^2\) solution always exist?}
For a classical solution to exist, the boundary condition \( f(\theta) \) must be sufficiently smooth to ensure uniform convergence of its Fourier expansion. Specifically, \( f(\theta) \) should be at least twice differentiable. When \( f \) is \( C^2 \), the Fourier expansion converges uniformly, which guarantees the existence of a well-behaved solution inside. This smoothness ensures that the solution satisfies the Laplace equation, which describes harmonic functions. Even if \( f \) is not twice differentiable, a solution can still exist, but with lower regularity. Also, the solution inside the disc is regularized by noise because the Laplace operator is closely tied to Brownian motion. As a result, even if the boundary function \( f \) has some roughness, the solution \( u \) inside the disc can still be differentiated as many times as necessary. However, if \( f \) is not smooth, then the solution near the boundary may reflect that, hence limiting the regularity to the level of \( f \).  

\subsubsection*{The Role of the Disc}
The disc geometry makes sure that the boundary is smooth and compact, which allows for the separation of variables and Fourier expansion; this simplifies solving the equation and ensures the uniqueness of the solution! For example, Brownian motion cannot fully explore regions with spikes or discontinuities because it follows random, continuous paths that are unlikely to precisely hit sharp or narrow areas. In smooth domains, Brownian motion naturally diffuses across the entire boundary; meanwhile when the boundary has a spike or discontinuity, the probability of the random path intersecting or reaching these points is significantly lower resulting in less regular solutions to differential equations.

\subsubsection{The Angular Equation}

First, we solve the angular equation
\[
\frac{d^2 \Theta}{d\theta^2} + n^2 \Theta = 0
\]
This is a standard linear second-order ordinary differential equation (ODE) with constant coefficients and its general solution is
\[
\Theta(\theta) = A_n \cos(n\theta) + B_n \sin(n\theta)
\]
where \(A_n\) and \(B_n\) are constants. An important consideration is that \(\Theta(\theta)\) must be \(2\pi\)-periodic because \(\theta\) represents the angular coordinate in polar coordinates, which is periodic with period \(2\pi\). That is, \(\Theta(\theta) = \Theta(\theta + 2\pi)\). This periodicity condition imposes restrictions on \(n\).

\[
\Theta(\theta + 2\pi) = A_n \cos(n(\theta + 2\pi)) + B_n \sin(n(\theta + 2\pi)) = \Theta(\theta)
\]
\[
\cos(n(\theta + 2\pi)) = \cos(n\theta + 2\pi n) = \cos(n\theta)
\]
\[
\sin(n(\theta + 2\pi)) = \sin(n\theta + 2\pi n) = \sin(n\theta)
\]
However, these equalities hold if \(n\) is an integer, because \(\cos(2\pi n) = \cos(0) = 1\) and \(\sin(2\pi n) = \sin(0) = 0\) when \(n\) is an integer. Therefore, to satisfy the periodicity condition, \(n\) must be an integer. This shows that \(n^2\) is the square of a natural number.

\subsubsection{The Radial Equation}

Next, we solve the radial equation
\[
r^2 \frac{d^2 R}{dr^2} + r \frac{dR}{dr} - n^2 R = 0
\]
This is a Cauchy-Euler equation, which can be solved using the trial solution \(R(r) = r^m\), where \(m\) is a constant to be determined. Substituting \(R(r) = r^m\) into the radial equation
\[
r^2 \left( m(m-1) r^{m-2} \right) + r \left( m r^{m-1} \right) - n^2 r^{m} = 0
\]
\[
m(m-1) r^{m} + m r^{m} - n^2 r^{m} = 0
\]
\[
\left[ m(m-1) + m - n^2 \right] r^{m} = 0
\]
\[
\left[ m^2 - n^2 \right] r^{m} = 0
\]
For a nontrivial solution, we need
\[
m^2 - n^2 = 0 \implies m = \pm n
\]
Therefore, the general solution for \(R(r)\) is
\[
R(r) = C_n r^n + D_n r^{-n}
\]
where \(C_n\) and \(D_n\) are constants.

\subsection{Combining Solutions and Applying Boundary Conditions}

Combining the solutions for \(\Theta(\theta)\) and \(R(r)\), the general solution to the Laplace equation in the disc is
\[
u(r, \theta) = \sum_{n=0}^{\infty} \left( C_n r^n + D_n r^{-n} \right) \left( A_n \cos(n\theta) + B_n \sin(n\theta) \right)
\]
Since we are solving inside a disc (i.e., \(0 \leq r \leq R\)), and to ensure that \(u(r, \theta)\) remains finite at \(r = 0\), we must remove the \(r^{-n}\) terms for \(n > 0\), as they become singular at \(r = 0\). So, the solution simplifies
\[
u(r, \theta) = C_0 + \sum_{n=1}^{\infty} r^n \left( A_n \cos(n\theta) + B_n \sin(n\theta) \right)
\]
Now, we apply the boundary condition \(u(R, \theta) = f(\theta)\). Setting \(r = R\), we get
\[
u(R, \theta) = C_0 + \sum_{n=1}^{\infty} R^n \left( A_n \cos(n\theta) + B_n \sin(n\theta) \right) = f(\theta)
\]
This is a Fourier series expansion of \(f(\theta)\) in terms of sines and cosines. The coefficients \(C_0\), \(A_n\), and \(B_n\) can be determined by expanding \(f(\theta)\) in a Fourier series
\[
f(\theta) = a_0 + \sum_{n=1}^{\infty} \left( a_n \cos(n\theta) + b_n \sin(n\theta) \right)
\]
By matching coefficients, we have
\[
C_0 = a_0, \quad A_n = \frac{a_n}{R^n}, \quad B_n = \frac{b_n}{R^n}
\]

\subsection{Uniqueness of Solution}

An important question is whether the solution to the Laplace equation with given boundary conditions is unique. 

\subsubsection{Proof of the Maximum Principle}

Let \( w(r, \theta) \) be a harmonic function defined on a closed domain \( \overline{D} \) (including its boundary \( \partial D \)). If \( w \) attains its maximum or minimum value at an interior point of \( D \), then \( w \) is constant throughout \( D \).

\vspace{2mm}
\textbf{Case 1: Maximum in the Interior}

Suppose \( w \) is harmonic in \( D \) (i.e., \( \nabla^2 w = 0 \)) and reaches a local maximum at some interior point \( (r_0, \theta_0) \). Consider the function \( v(r, \theta) = w(r, \theta) - \epsilon r^2 \) for a small \(\epsilon < 0\). Assume \( v(r, \theta) \) attains a local maximum at the same point \( (r_0, \theta_0) \). At this maximum point, the first-order partial derivatives of \( v \) disappear
\[
\frac{\partial v}{\partial r}(r_0, \theta_0) = \frac{\partial w}{\partial r}(r_0, \theta_0) - 2\epsilon r_0 = 0
\]
\[
\frac{\partial v}{\partial \theta}(r_0, \theta_0) = \frac{\partial w}{\partial \theta}(r_0, \theta_0) = 0
\]

The second-order partial derivatives satisfy
\[
\frac{\partial^2 v}{\partial r^2}(r_0, \theta_0) = \frac{\partial^2 w}{\partial r^2}(r_0, \theta_0) - 2\epsilon
\]
\[
\frac{\partial^2 v}{\partial \theta^2}(r_0, \theta_0) = \frac{\partial^2 w}{\partial \theta^2}(r_0, \theta_0)
\]

The Laplacian of \( v \) at \( (r_0, \theta_0) \) then is
\[
\nabla^2 v(r_0, \theta_0) = \frac{\partial^2 v}{\partial r^2} + \frac{1}{r_0} \frac{\partial v}{\partial r} + \frac{1}{r_0^2} \frac{\partial^2 v}{\partial \theta^2}
\]

Using the conditions at the maximum point
\[
\nabla^2 v(r_0, \theta_0) = \frac{\partial^2 w}{\partial r^2} - 2\epsilon + \frac{1}{r_0^2} \frac{\partial^2 w}{\partial \theta^2}
\]

Since \( \nabla^2 v = 4\epsilon \) and \( \epsilon < 0 \), we have
\[
\frac{\partial^2 w}{\partial r^2}(r_0, \theta_0) - 2\epsilon + \frac{1}{r_0^2} \frac{\partial^2 w}{\partial \theta^2} = 4\epsilon
\]

It follows that
\[
\frac{\partial^2 w}{\partial r^2}(r_0, \theta_0) + \frac{1}{r_0^2} \frac{\partial^2 w}{\partial \theta^2} = 2\epsilon
\]

However, because \( w \) is harmonic (\( \nabla^2 w = 0 \)), the left-hand side must equal zero. This shows that \( w \) cannot attain a maximum at an interior point unless \( w \) is constant.

\vspace{2mm}
\textbf{Case 2: Minimum in the Interior}

Assume \( w(r, \theta) \) attains a local minimum at an interior point \( (r_0, \theta_0) \). Consider the function \( v(r, \theta) = w(r, \theta) + \epsilon r^2 \), where \( \epsilon < 0 \). By similar reasoning, \( v \) is harmonic, but:
\[
\nabla^2 v = -4\epsilon < 0
\]

At the minimum point \( (r_0, \theta_0) \), the first-order partial derivatives again disappear
\[
\frac{\partial v}{\partial r}(r_0, \theta_0) = \frac{\partial w}{\partial r}(r_0, \theta_0) + 2\epsilon r_0 = 0
\]
\[
\frac{\partial v}{\partial \theta}(r_0, \theta_0) = \frac{\partial w}{\partial \theta}(r_0, \theta_0) = 0
\]

The second-order partial derivatives satisfy
\[
\frac{\partial^2 v}{\partial r^2}(r_0, \theta_0) = \frac{\partial^2 w}{\partial r^2}(r_0, \theta_0) + 2\epsilon
\]
\[
\frac{\partial^2 v}{\partial \theta^2}(r_0, \theta_0) = \frac{\partial^2 w}{\partial \theta^2}(r_0, \theta_0)
\]

The Laplacian of \( v \) at \( (r_0, \theta_0) \) is
\[
\nabla^2 v(r_0, \theta_0) = \frac{\partial^2 w}{\partial r^2} + 2\epsilon + \frac{1}{r_0^2} \frac{\partial^2 w}{\partial \theta^2}
\]

Since \( \nabla^2 v = -4\epsilon \) and \( \epsilon < 0 \), we have
\[
\frac{\partial^2 w}{\partial r^2}(r_0, \theta_0) + \frac{1}{r_0^2} \frac{\partial^2 w}{\partial \theta^2} = 2\epsilon
\]

This is a contradiction, as \( w \) being harmonic requires \( \nabla^2 w = 0 \). So, \( w \) cannot attain a minimum at an interior point unless \( w \) is constant. Therefore, \( w \) attains its maximum or minimum at an interior point of \( D \) as \( w \) must be constant throughout \( D \); so, the extreme values of a non-constant harmonic function on \( D \) occur only on the boundary \( \partial D \). 

\subsubsection{Showing \( w \equiv 0 \)}

Suppose \( w \) is harmonic in \( D \) and satisfies the boundary condition \( w = 0 \) on \( \partial D \). By the maximum principle, the maximum and minimum of \( w \) occur on \( \partial D \). Since \( w = 0 \) on \( \partial D \), the maximum and minimum values of \( w \) are both zero. Therefore, \( w \) cannot attain any value other than zero within \( D \) unless it is constant. We know that \( w \) is harmonic and the maximum and minimum values are the same, so it follows that \( w \) is constant throughout \( D \). Given \( w = 0 \) on \( \partial D \), 
\[
w(r, \theta) \equiv 0 \quad \text{in} \quad D
\]

\subsubsection{Proof of Uniqueness}

Suppose there are two solutions \( u_1(r,\theta) \) and \( u_2(r,\theta) \) to the Laplace equation \( \nabla^2 u = 0 \) in the disc \( D \) of radius \( R \), both satisfying boundary condition \( u(R,\theta) = f(\theta) \). Let \( w(r,\theta) = u_1(r,\theta) - u_2(r,\theta) \) s.t. \( w(r,\theta) \) satisfies

\begin{enumerate}
    \item Laplace equation: \( \nabla^2 w = \nabla^2 u_1 - \nabla^2 u_2 = 0 - 0 = 0 \)
    \item Boundary condition: \( w(R,\theta) = u_1(R,\theta) - u_2(R,\theta) = f(\theta) - f(\theta) = 0 \)
    
\end{enumerate} Therefore, \( w(r,\theta) \) is a harmonic function in \( D \) that disappears on boundary \( r = R \). By the maximum principle for harmonic functions, we have that the maximum of \( w \) in the closed domain \( \overline{D} \) is achieved on the boundary \( \partial D \) and, as \( w \) is continuous on \( \overline{D} \) and harmonic in \( D \) and \( w \) vanishes on the boundary, then it follows that the maximum and minimum of \( w \) are zero. So, \( w(r,\theta) = 0 \) \(\forall\) \( (r,\theta) \in D \) and, therefore, \( u_1(r,\theta) = u_2(r,\theta) \) throughout \( D \), proving that the solution is unique.

\subsection{Alternative Methods}

If we don't use separation of variables, an alternative is the Poisson integral formula.

\[
u(r, \theta) = \frac{1}{2\pi} \int_{0}^{2\pi} \frac{R^2 - r^2}{R^2 - 2rR \cos(\theta - \phi) + r^2} f(\phi) \, d\phi
\] To find the Poisson integral, we solve the Laplace equation \( \nabla^2 u = 0 \) in polar coordinates under the Dirichlet boundary condition \( u(R, \theta) = f(\theta) \). The approach uses the Green's function for the Laplacian in a disc, which is defined to satisfy the Laplace equation everywhere except at a singular point. The Green's function for a point located at \( (r, \theta) \) in a disc of radius \( R \) is
   \[
   G(r, \theta; r', \phi) = -\frac{1}{2\pi} \ln \left( \frac{1 - 2rr'\cos(\theta - \phi) + r^2r'^2}{r^2 + r'^2 - 2rr'\cos(\theta - \phi)} \right)
   \]
   Applying boundary condition \( G = 0 \) on \( r = R \) makes sure that the solution satisfies the Dirichlet problem. For a harmonic function \( u \), applying Green's identity with the Green's function gives
   \[
   u(r, \theta) = \frac{1}{2\pi} \int_{0}^{2\pi} G(R, \phi; r, \theta) f(\phi) \, d\phi
   \]

\subsubsection{Proof of Equivalence between the Poisson Integral and Separation of Variables}

Now we can show the Poisson integral formula gives the same solution for \( u(r, \theta) \) as the separation of variables method.

\[
u(r, \theta) = \frac{1}{2\pi} \int_{0}^{2\pi} \frac{R^2 - r^2}{R^2 - 2 r R \cos(\theta - \phi) + r^2} f(\phi) \, d\phi
\]

Expanding the denominator using the identity

\[
\frac{1}{R^2 - 2 r R \cos(\theta - \phi) + r^2} = \sum_{n=0}^{\infty} \left( \frac{r}{R} \right)^n \cos(n(\theta - \phi))
\]

Substituting this into the integral

\[
u(r, \theta) = \frac{1}{2\pi} \int_{0}^{2\pi} \left( 1 + 2 \sum_{n=1}^{\infty} \left( \frac{r}{R} \right)^n \cos(n(\theta - \phi)) \right) f(\phi) \, d\phi
\]

Applying the Fourier series representation of \( f(\phi) \)

\[
f(\phi) = \frac{a_0}{2} + \sum_{n=1}^{\infty} \left( a_n \cos(n\phi) + b_n \sin(n\phi) \right)
\]

We have 

\[
u(r, \theta) = \frac{a_0}{2} + \sum_{n=1}^{\infty} \left( \frac{r}{R} \right)^n \left( a_n \cos(n\theta) + b_n \sin(n\theta) \right)
\]

This is the same solution as through separation of variables, which confirms that both methods give us the same harmonic function inside the disc; this follows from the uniqueness of the solution to Laplace’s equation with the same boundary conditions.

\subsection{Example Problem}

Suppose we are given the boundary condition
\[
u(R, \theta) = f(\theta) = \sin(2\theta)
\]
where \( R \) is the radius of the disc. Our goal is to find \( u(r, \theta) \) inside the disc that satisfies the Laplace equation \( \nabla^2 u = 0 \) and the given boundary condition. From the general solution
\[
u(r, \theta) = C_0 + \sum_{n=1}^{\infty} r^n \left( A_n \cos(n\theta) + B_n \sin(n\theta) \right)
\]

Since \( f(\theta) = \sin(2\theta) \), the only nonzero Fourier coefficient is \( b_2 = 1 \).

Matching the boundary condition at \( r = R \)
\[
u(R, \theta) = C_0 + R^2 B_2 \sin(2\theta) = \sin(2\theta)
\]
This implies that \( C_0 = 0 \) and \( B_2 = \dfrac{1}{R^2} \). So, the solution inside the disc is
\[
u(r, \theta) = r^2 B_2 \sin(2\theta) = \left( \dfrac{r^2}{R^2} \right) \sin(2\theta)
\]

This function satisfies the Laplace equation and the boundary condition \( u(R, \theta) = \sin(2\theta) \).

\section{Holomorphic Functions and the Laplace Equation}

Holomorphic functions always satisfy the Laplace equation, but not all solutions to the Laplace equation are holomorphic. To demonstrate this, consider a holomorphic function \(f(z) = u(x, y) + iv(x, y)\), where \(z = x + iy\) and \(u(x, y)\) and \(v(x, y)\) are the real and imaginary parts of the function. For \(f(z)\) to be holomorphic, the real and imaginary components must satisfy the Cauchy-Riemann equations
\[
\frac{\partial u}{\partial x} = \frac{\partial v}{\partial y}, \quad \frac{\partial u}{\partial y} = -\frac{\partial v}{\partial x}
\]
We now show that both \(u(x, y)\) and \(v(x, y)\) satisfy the Laplace equation \(\nabla^2 u = 0\). Differentiating the Cauchy-Riemann equation \(\frac{\partial u}{\partial x} = \frac{\partial v}{\partial y}\) with respect to \(x\) gives \(\frac{\partial^2 u}{\partial x^2} = \frac{\partial^2 v}{\partial y \partial x}\). Similarly, differentiating \(\frac{\partial u}{\partial y} = -\frac{\partial v}{\partial x}\) with respect to \(y\) gives us \(\frac{\partial^2 u}{\partial y^2} = -\frac{\partial^2 v}{\partial x \partial y}\). Adding these two equations gives us \(\nabla^2 u = \frac{\partial^2 u}{\partial x^2} + \frac{\partial^2 u}{\partial y^2} = 0\). Hence, the real and imaginary parts of any holomorphic function satisfy the Laplace equation.

However, the reverse is not true as not all solutions to the Laplace equation correspond to holomorphic functions. A real-valued function that satisfies \(\nabla^2 u = 0\) is called harmonic, but it does not necessarily have an associated holomorphic function. For example, the function \(u(x, y) = x^2 - y^2\) is harmonic, since \(\nabla^2 u = \frac{\partial^2 u}{\partial x^2} + \frac{\partial^2 u}{\partial y^2} = 2 - 2 = 0\), but it does not correspond to the real part of a holomorphic function. There is no function \(v(x, y)\) that can satisfy the Cauchy-Riemann equations alongside \(u(x, y)\). This shows that while holomorphic functions must satisfy the Laplace equation, not all solutions to the Laplace equation are holomorphic.

\subsection{Mean Value Property}

Let \( u \) be a harmonic function in an open domain \( \Omega \subset \mathbb{R}^n \). For any closed ball \( \overline{B}(x_0, r) \subset \Omega \), the property states
\[
u(x_0) = \frac{1}{|\partial B_r|} \int_{\partial B(x_0, r)} u(y) \, d\sigma(y)\]
where \( \partial B(x_0, r) \) is the sphere of radius \( r \) centered at \( x_0 \), \( |\partial B_r| \) is the surface area of the sphere, and \( d\sigma(y) \) is the surface measure on \( \partial B(x_0, r) \). Since \( u \) is harmonic in \( \Omega \), it satisfies Laplace's equation
\[
\nabla^2 u = 0 \quad \text{in } B(x_0, r)
\]
We consider the integral of \( u \) over the surface \( \partial B(x_0, r) \)
\[
A(r) = \frac{1}{|\partial B_r|} \int_{\partial B(x_0, r)} u(y) \, d\sigma(y)
\]

By symmetry, \( A(r) \) depends only on the value of \( u \) at points \( y \) on the sphere and their radial distance from \( x_0 \). We write \( u(x) \) in terms of its integral representation using the Poisson kernel for the ball such that we have
\[
u(r, \theta) = \frac{1}{|\partial B_r|} \int_{\partial B(x_0, r)} P(r, \theta - \phi) u(\phi) \, d\sigma(\phi)\]
where the Poisson kernel in polar coordinates is
\[
P(r, \theta - \phi) = \frac{R^2 - r^2}{R^2 - 2rR\cos(\theta - \phi) + r^2}
\]
The Poisson kernel has rotational symmetry as it depends only on \( |\theta - \phi| \), the angle between \( x \) and \( y \) on the sphere. Expanding \( P(r, \theta - \phi) \) in terms of its Fourier series, we have
\[
P(r, \theta - \phi) = 1 + 2 \sum_{n=1}^\infty \left( \frac{r}{R} \right)^n \cos(n(\theta - \phi))
\]

The integral of \( P(r, \theta - \phi) \) over the sphere simplifies because of orthogonality of the trigonometric terms; so, the constant term \( 1 \) is the only term remaining
\[
\int_{\partial B(x_0, r)} P(r, \theta - \phi) \, d\sigma(\phi) = |\partial B_r|
\]

So, \( P(r, \theta - \phi) \) integrates to 1 over the sphere, showing that the mean value property holds. 
\[
\phi'(r) = \frac{\partial}{\partial r} \left( \frac{1}{|\partial B_r|} \int_{\partial B(x_0, r)} u(y) \, d\sigma(y) \right) = 0 
\]
This shows \( u \) is constant along radial directions on \( \partial B(x_0, r) \). Using symmetry of \( P(r, \theta - \phi) \), the value of \( u(x) \) at the center \( x_0 \) is
\[
u(x_0) = \frac{1}{|\partial B_r|} \int_{\partial B(x_0, r)} u(y) \, d\sigma(y)\]
where \( u(y) \) depends only on \( r \). This shows that \( u(x_0) \) is equal to the average value of \( u \) over the sphere. So, the symmetry of the Poisson kernel and \( u \) being harmonic ensures that the mean value property holds. So, we have
\[
u(x_0) = \frac{1}{|\partial B_r|} \int_{\partial B(x_0, r)} u(y) \, d\sigma(y)\]

\subsubsection{The Value is Constant}
For the Fourier series expansion
\[
P(r, \theta - \phi) = 1 + 2 \sum_{n=1}^\infty \left( \frac{r}{R} \right)^n \cos(n(\theta - \phi))
\]
only the constant term \( 1 \) is left after integration because
\[
\int_{\partial B(x_0, r)} \cos(n(\theta - \phi)) \, d\sigma(\phi) = 0 \quad \text{for } n \geq 1
\]
So, we have
\[
\int_{\partial B(x_0, r)} P(r, \theta - \phi) \, d\sigma(\phi) = |\partial B_r|
\]

Substituting this back into the integral for \( u(x_0) \)
\[
u(x_0) = \frac{1}{|\partial B_r|} \int_{\partial B(x_0, r)} u(y) \, d\sigma(y)
\]

This shows \( P(r, \theta - \phi) \) integrates to \( 1 \) over the sphere and, therefore, that the mean value property holds. Using symmetry of \( P(r, \theta - \phi) \), we have the value of \( u(x) \) at the center \( x_0 \) as
\[
u(x_0) = \frac{1}{|\partial B_r|} \int_{\partial B(x_0, r)} u(y) \, d\sigma(y)
\]
where \( u(y) \) depends only on \( r \). This shows that the value is constant for the mean value theorem as we have

\[
u(x_0) = \frac{1}{|\partial B_r|} \int_{\partial B(x_0, r)} u(y) \, d\sigma(y)
\]

\subsubsection{Mean Value for Harmonic Functions}

Harmonic functions, solutions to the Laplace equation \( \nabla^2 u = 0 \), satisfy the mean value property.

\vspace{2mm}
Assume \( u \) is harmonic. Since \( u \) satisfies \( \nabla^2 u = 0 \), it is twice continuously differentiable in \( D \).
   Set \( x = x_0 + r\cos\theta \) and \( y = y_0 + r\sin\theta \). The Laplace equation becomes
   \[
   \frac{\partial^2 u}{\partial r^2} + \frac{1}{r} \frac{\partial u}{\partial r} + \frac{1}{r^2} \frac{\partial^2 u}{\partial \theta^2} = 0
   \]
   Consider average value of \( u \) on the boundary of the disc \( r = R \)
   \[
   \frac{1}{2\pi} \int_{0}^{2\pi} u(x_0 + R\cos\theta, y_0 + R\sin\theta) \, d\theta
   \]
   By applying Green’s identities and the harmonicity of \( u \), the radial average of \( u \) over concentric circles is constant
   \[
   \frac{1}{2\pi} \int_{0}^{2\pi} u(x_0 + r\cos\theta, y_0 + r\sin\theta) \, d\theta = u(x_0, y_0)
   \]
   Because \( u \) is harmonic, any deviation from the average in the angular direction contributes zero net effect when integrated over a full circle, ensuring that the value at the center equals the average over the boundary. So, \( u(x_0, y_0) \) equals the average of \( u \) over the boundary of the disc. Going the other way, assume \( u \) satisfies the mean value property; then for any disc \( D \) of radius \( R \)
\[
u(x_0, y_0) = \frac{1}{2\pi} \int_{0}^{2\pi} u(x_0 + R\cos\theta, y_0 + R\sin\theta) \, d\theta
\]
Differentiating this expression under the integral sign and applying properties of harmonic functions shows \( \nabla^2 u = 0 \). Therefore, \( u \) is harmonic.

\subsection{Converse to Mean Value Theorem}

If \( u \in C^2(U) \) satisfies
\[
u(x) = \fint_{\partial B(x,r)} u \, dS
\]
for each ball \( B(x,r) \subset U \), then \( u \) is harmonic.

\begin{proof}
If \( \Delta u \not\equiv 0 \), there exists some ball \( B(x,r) \subset U \) such that, say, \( \Delta u > 0 \) within \( B(x,r) \). Consider a radial function \( \phi \) as before. Then,
\[
0 = \phi'(r) = \frac{1}{|B(x,r)|} \int_{B(x,r)} \Delta u(y) \, dy
\]
which leads to a contradiction if \( \Delta u(y) > 0 \). Therefore, \( \Delta u = 0 \) in \( U \), and \( u \) is harmonic.
\end{proof}

\subsection{Strong Maximum Principle and Uniqueness}

Suppose \( u \in C^2(U) \cap C(\overline{U}) \) is harmonic within \( U \). Then
\begin{enumerate}
    \item[(i)] \[
    \sup_{U} u = \sup_{\partial U} u
    \]
    \item[(ii)] If \( U \) is connected and there exists a point \( x_0 \in U \) such that
    \[
    u(x_0) = \sup_{\overline{U}} u
    \]
    then \( u \) is constant within \( U \).
\end{enumerate}

\begin{proof}
Suppose that there exists a point \( x_0 \in U \) such that \( u(x_0) = M := \sup_{U} u \). For \( 0 < r < \mathrm{dist}(x_0, \partial U) \), the mean value property tells us that
\[
M = u(x_0) = \fint_{\partial B(x_0,r)} u \, dS \leq M
\]

As equality holds only if \( u \equiv M \) in \( B(x_0, r) \), we have \( u(y) = M \) for all \( y \in B(x_0, r) \). Thus, the set \( \{ x \in U \mid u(x) = M \} \) is open and relatively closed in \( U \), and therefore equals \( U \) if \( U \) is connected.
\end{proof}

\begin{remark}
The strong maximum principle asserts in particular that if \( U \) is connected and \( u \in C^2(U) \cap C(\overline{U}) \) satisfies
\[
\begin{cases}
\Delta u = 0 & \text{in } U \\
u = g & \text{on } \partial U
\end{cases}
\]
where \( g \geq 0 \), then \( u \) is positive \emph{everywhere} in \( U \) if \( g \) is positive \emph{somewhere} on \( \partial U \).
\end{remark}

\section{Definition and Properties of Brownian Motion}

According to Mörters and Perez, Brownian motion is a stochastic process that models continuous random motion, denoted \( (B_t)_{t \geq 0} \). It has the following properties

\vspace{2mm}
\begin{itemize}
    \item \textbf{Initial Condition:} \( B_0 = 0 \) makes sure that the motion starts at the origin.

    \item \textbf{Independence:} For \( 0 \leq t_1 < t_2 < \dots < t_n \), the increments \( B_{t_2} - B_{t_1}, B_{t_3} - B_{t_2}, \dots, B_{t_n} - B_{t_{n-1}} \) are independent, reflecting how the future evolution of the process does not depend on its past behavior.

    \item \textbf{Stationarity:} The distribution of increments \( B_{t+s} - B_t \) depends only on the time difference \( s \), not on \( t \) itself, which makes sure the process is time homogeneous.

    \item \textbf{Gaussian Increments:} For \( s > 0 \), increments \( B_{t+s} - B_t \) are normally distributed with mean \( 0 \) and variance \( s \): \( B_{t+s} - B_t \sim \mathcal{N}(0, s) \). This shows the structure of Brownian motion as a Gaussian process.

    \item \textbf{Continuity of Paths:} The function \( t \mapsto B_t \) is continuous; however, while the paths are continuous, they are also highly irregular.
\end{itemize}

\subsection{Stopping Time}
In stochastic processes, a stopping time \( \tau \) is a random variable that represents the time at which a condition is met (determined by the history of the process up to time \( \tau \)). Specifically, \( \tau \) is a stop time with respect to a filtration \( (\mathcal{F}_t)_{t \geq 0} \) if, for every \( t \geq 0 \), the event \( \{\tau \leq t\} \) belongs to \( \mathcal{F}_t \), making sure that the decision to stop can be made based only on the information available up to \( t \). Stopping times are important in applications involving Brownian motion as the stopping time corresponds to the random moment when a Brownian path first hits a specified boundary.

\subsection{Value on the Boundary Using Brownian Motion}
To find the value of a harmonic function \( u \) at a point \( x_0 \) in a domain \( \Omega \), we can simulate a Brownian motion starting from \( x_0 \). The procedure is 

\begin{itemize}
\item[1.] Start a Brownian motion \( B_t \) at \( x_0 \) and let it randomly evolve 
\item[2.] Stop the motion for the first time \( \tau \) when it hits boundary \( \partial \Omega \)
\item[3.] Record the value of \( u \) at the point where the Brownian motion hits \( u(B_\tau) \)
\end{itemize}

The mean value of \( u(B_\tau) \), averaged over all possible paths, gives
\[
u(x_0) = \mathbb{E}[u(B_\tau)]
\]
where \( \tau \) is the stopping time when \( B_t \) hits the boundary \( \partial \Omega \). 

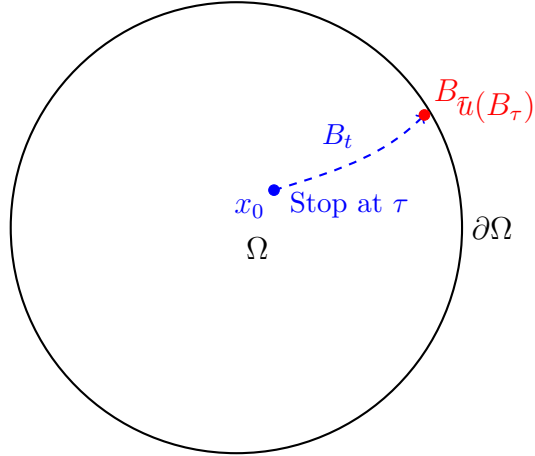
\vspace{2mm}
\begin{center}
\begin{tikzpicture}[scale=1]

\draw[thick] (0,0) circle (3cm) node[below right, font=\Large] {$\Omega$};

\draw[thick] (3,0) node[right, font=\Large] {$\partial \Omega$};

\filldraw[blue] (0.5,0.5) circle (2pt) node[below left, font=\Large] {$x_0$};

\draw[thick, dashed, ->, blue] (0.5,0.5) 
  .. controls (1.5,0.8) and (2,1) .. (2.5,1.5) 
  node[midway, above left, font=\large, blue] {$B_t$};

\filldraw[red] (2.5,1.5) circle (2pt) node[above right, font=\Large, red] {$B_\tau$};

\node[red, right, font=\Large] at (2.8,1.6) {$u(B_\tau)$};

\node[blue, below, font=\large] at (1.5,0.6) {Stop at $\tau$};

\end{tikzpicture}
\end{center}

\subsection{Solving the Laplace Equation with Brownian Motion}

For any domain \( \Omega \subset \mathbb{R}^n \), the Laplace equation is 
\[
\nabla^2 u = 0 \quad \text{in } \Omega
\]
with boundary condition \( u = f \) on \( \partial \Omega \). According to Morters and Perez' book, the solution is
\[
u(x) = \mathbb{E}[f(B_\tau)]
\]
where \( B_t \) is a Brownian motion starting at \( x \), and \( \tau \) is the stopping time when \( B_t \) hits the boundary \( \partial \Omega \). For a disc \( \Omega = \{ x \in \mathbb{R}^2 : |x| < R \} \), the solution simplifies because of the symmetry of the domain. The value \( u(x) \) depends only on the boundary values and the geometry of the disc.

\subsubsection{Example: Finding Temperature Inside a Disc}

Suppose we have disc of radius \( R = 1 \). The temperature at the boundary is
\[
f(\theta) = 100 \sin(2\theta)
\] and we want to calculate the temperature \( u(r, \theta) \) at any point \( (r, \theta) \) inside the disc using Brownian motion. We know that the temperature distribution inside the disc satisfies the Laplace equation
\[
\nabla^2 u = 0
\]  
This implies \( u(r, \theta) \) is harmonic inside the disc. Imagine heat particles diffusing randomly inside the disc; a Brownian particle starts at \( (r, \theta) \) and moves randomly until it hits the boundary at \( (1, \phi) \). The temperature at the point \( (r, \theta) \) is the average of the boundary temperatures the particle encounters during this process. This is expressed as

\[
u(r, \theta) = \mathbb{E}[f(B_{\tau})]
\]  
where \( B_{\tau} \) is the location on the boundary where the particle first exits, and \( \tau \) is the stopping time when the particle reaches the boundary. Though we would need software to compute these values, an outline of the process follows. First, choose a point inside disc \( (r, \theta) \), suppose \( r = 0.5 \) and \( \theta = \frac{\pi}{4} \). From \( (0.5, \frac{\pi}{4}) \), generate a random Brownian path; the particle will move in small random steps until it exits the disc at some boundary point \( (1, \phi) \). We will record the angle \( \phi \) where the particle exits. At the boundary, we are given that the temperature is given by \( f(\phi) = 100 \sin(2\phi) \). So, for each Brownian path, we can calculate the temperature at the exit point \( f(\phi) \) and average the values of \( f(\phi) \) obtained from multiple trials of Brownian paths. Suppose we simulate three Brownian paths from \( (0.5, \frac{\pi}{4}) \).

\[
\begin{array}{|c|c|c|}
\hline
\text{Path \#} & \text{Boundary Exit Angle } \phi & \text{Boundary Temperature } f(\phi) \\
\hline
1 & \frac{\pi}{3} & 100 \sin\left(2 \cdot \frac{\pi}{3} \right) = 87 \\
\hline
2 & \frac{5\pi}{6} & 100 \sin\left(2 \cdot \frac{5\pi}{6} \right) = 50 \\
\hline
3 & \frac{7\pi}{4} & 100 \sin\left(2 \cdot \frac{7\pi}{4} \right) = -70 \\
\hline
\end{array}
\]

\[
u(0.5, \frac{\pi}{4}) \approx \frac{87 + 50 + (-70)}{3} = 22.33
\]

By repeating this for thousands of paths (hence why this would need to be done with a software), the average converges to the true temperature at \( (0.5, \frac{\pi}{4}) \).

\subsection*{Remark}
To compute the expected value of a harmonic function, we need a probability density that describes how likely it is for a point to be hit at the boundary. The Poisson kernel gives us this density as it represents the solution to the Dirichlet problem by describing how values inside a domain are influenced by values on the boundary. The Poisson kernel models the distribution of boundary points a Brownian particle is likely to hit when it first exits the domain. So, when \( r = 0 \) (at the center of the domain), the Poisson kernel simplifies to 1; this tells us that the probability of hitting any point on the boundary has a uniform distribution and, therefore, from the center, all boundary points are equally as likely to be hit.  

\section{Applications of the Laplace Equation in Imaging}

\subsection*{What is Imaging?}
Imaging refers to techniques that create visual representations of physical phenomena, allowing interpretation, analysis, or diagnostics. These methods rely on data capture, processing, and visualization to produce visual outputs.

\subsection*{Why is the Laplace Equation Useful in Imaging?}
The Laplace equation is important as it can ensure smoothness and consistency in many imaging tasks. The Laplace equation is often most useful in the processing stage, as it provides a framework for interpolating missing data, reducing noise, and enhancing resolution across various domains.

\begin{enumerate}
    \item \textbf{Noise Reduction:} Removing irregularities while preserving critical features.
    \item \textbf{Data Interpolation:} Filling in missing or corrupt regions.
    \item \textbf{Enhancing Resolution:} Smoothing transitions to improve clarity.
\end{enumerate}

\subsubsection*{Step 1: Noise reduction} Noise is inherent in most imaging systems and can obscure important details. The Laplace equation smooths these irregularities by enforcing harmonicity, ensuring that each point in the image is consistent with its surroundings.
\vspace{1mm}

\textbf{Outline of a Proof:}
Let \( u(x, y) \) represent the intensity of the image. If we can minimize this equation
\[
E(u) = \int_\Omega \left| \nabla u \right|^2 \, d\Omega
\]
to get the Euler-Lagrange equation \( \nabla^2 u = 0 \), the solution of this equation guarantees that the intensity at each point is the average of its neighbors, which reduces noise.
\vspace{1mm}

\textbf{Reference:} Perona and Malik, \emph{Scale-Space and Edge Detection Using Anisotropic Diffusion}, IEEE Transactions on Pattern Analysis, 1990.
\vspace{2mm}

\subsubsection*{Step 2: Data Interpolation}
In imaging, incomplete or corrupted data needs to be reconstructed. The Laplace equation gives a consistent method to fill these gaps.
\vspace{1mm}

\textbf{Outline of a Proof:}
Let \( u(x, y) \) be the image with missing data in a region \( D \) bounded by \( \partial D \). We can solve
\[
\nabla^2 u = 0 \quad \text{in } D, \quad u = g \quad \text{on } \partial D
\]
By the mean value property, we have
\[
u(x_0, y_0) = \frac{1}{2\pi r} \int_{\partial B((x_0, y_0), r)} u(x, y) \, dS
\]
where \( B((x_0, y_0), r) \) is a ball centered at \( (x_0, y_0) \). This ensures smooth transitions and fills in missing regions!
\vspace{1mm}

\textbf{Reference:} Haber and Oldenburg, \emph{Magnetic Resonance Image Reconstruction}, SIAM, 2007.
\vspace{2mm}

\subsubsection*{Step 3: Improving Resolution}
 Improving resolution means refining image details for better clarity; the Laplace equation smooths transitions while maintaining structural integrity.
\vspace{1mm}

\textbf{Outline of a Proof:}
Consider image intensity function \( u(x, y) \) with known boundary values \( g(x, y) \) on \( \partial \Omega \). We can solve
\[
\nabla^2 u = 0 \quad \text{in } \Omega, \quad u = g \quad \text{on } \partial \Omega
\]
Harmonic solutions ensure that the intensity values transition smoothly, creating a higher-resolution image.
\vspace{1mm}

\textbf{Reference:} Xu and Wang, \emph{Photoacoustic Imaging in Biomedicine}, CRC Press, 2010.
\vspace{2mm}

So, as seen, the Laplace equation is an important tool in imaging, as it guarantees smoothness and consistency. By addressing noise, missing data, and resolution challenges, it allows for robust imaging solutions.

\section{Example Imaging Problem Using the Laplace Equation}

Suppose we have a circular domain with radius $R$ representing a 2D image region and we know the boundary of this circle is  $f(\theta)$, $0 \le \theta < 2\pi$. Inside the circle, some portion of the image data is missing or needs to be smoothed. We want to find values such that each point inside the image be as \emph{harmonic} as possible, matching the boundary intensities at $r=R$ (a.k.a. we want to smooth the interior $\Omega$ using the Laplace equation). Below is a diagram showing the given circle.  

\begin{figure}[h!]
\centering
\begin{tikzpicture}[scale=1.3]
  \draw[thick] (0,0) circle (2cm);
  \filldraw (0,0) circle (1.2pt);
  \node[below right] at (0.2,-0.1) {$(0,0)$};
  \draw (0,0) -- (2,0);
  \node at (1,0.25) {$R$};
  \draw [->] (2.2,0.7) -- (1.6,1.3);
  \node at (2.6,0.7) {$u(R,\theta) = f(\theta)$};
  \draw[dashed] (-0.6,-0.1) .. controls (-0.4,0.5) and (0.4,0.5) .. (0.6,-0.1)
               .. controls (0.6,-0.4) and (-0.6,-0.5) .. (-0.6,-0.1);
  \node at (0,-0.8) {$\Omega$ (missing/corrupt region)};
\end{tikzpicture}
\caption{Circular imaging domain with boundary data $f(\theta)$ known at $r=R$. The dashed region $\Omega$ represents missing or corrupt data.}
\label{fig:imaging-problem}
\end{figure}
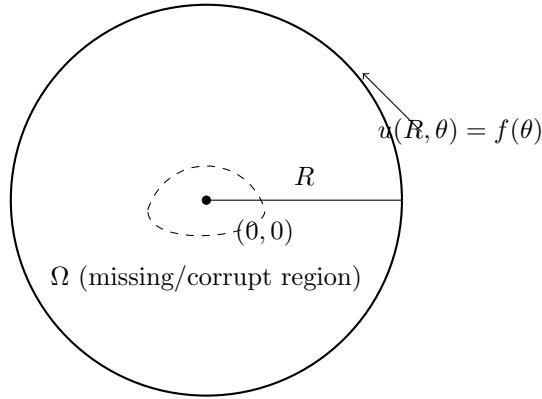

To smooth the image, let $u(r,\theta)$ be the image intensity at polar coordinates $(r,\theta)$, subject to the constraint that $u(R,\theta) = f(\theta)$. We have
\[
\nabla^2 u = 0,
\quad
\text{in the domain } 0 \le r < R
\quad
\text{with} 
\quad
u(R,\theta) = f(\theta)
\]

By theory explained earlier in this paper (separation of variables, Poisson integral formula, or Brownian motion arguments), we have the unique solution

\[
u(r,\theta) \;=\; \frac{1}{2\pi}
\int_{0}^{2\pi}
\frac{R^2 - r^2}{R^2 - 2Rr\cos(\theta - \phi) + r^2} \, f(\phi)\, d\phi
\]

This means $u(r,\theta)$ is smooth (a.k.a. harmonic) throughout the circle and exactly matches boundary data $f(\theta)$ at the outer boundary $r=R$. An alternative interpretation is to start a Brownian motion at the interior point $(r,\theta)$ and stop it when it first hits the boundary $r=R$. Then, the average boundary intensity observed over many random paths is exactly $u(r,\theta)$. This gives us a probabilistic method for solving the problem
\[
u(r,\theta)
\;=\;
\mathbb{E}\bigl[f\bigl(B_{\tau}\bigr)\bigr]
\]
where $\tau$ is the hitting time for the boundary of the disc.
\vspace{1mm}

\subsection{Remark}
Regularization by noise in the context of this imaging problem comes from the smoothing effect of the Laplace equation, driven by Brownian motion. When solving for the harmonic function \( u(r, \theta) \) inside the circular domain, the randomness of Brownian paths introduces a natural smoothing process. So, even if the boundary data \( f(\theta) \) is irregular or slightly discontinuous, the solution inside the disc is smooth because the noise inherent to Brownian motion averages out the inconsistencies, which then gives us a well-behaved harmonic function throughout the domain. This makes sure that  missing or corrupted data in the interior can be reconstructed smoothly, producing a clean and regularized image!

\section{Acknowledgments}
I would like to thank my mentor, Antonis Zitridis, for his invaluable knowledge, guidance, and support throughout this research project. I would also like to thank the University of Chicago Directed Reading Program Committee for providing me with the opportunity to produce and present this paper.

\end{document}